\documentclass[10pt]{article}
\begin{document}
\title{ {\bf Some open problems on cycles}
\thanks{ Revise title from 'Some unsolved problems on cycles'
To 'Some open problems on cycles' by the referees comments.
 Project supported by  the National Science Foundation of China (No. 11101358),
 NSF of Fujian(2011J01026), Innovation Program for Young Scientists of Fujian(2011J05014),
Fujian Provincial Training Foundation for "Bai-Quan-Wan Talents
Engineering",  IRTSTFJ, Project of Fujian Education Department (JA11165) and
Project of Zhangzhou Teachers College.}}
\author{{$Chunhui$ $Lai^{1,2}$, $Mingjing$ $Liu^{1}$}\\
{\small 1. Department of Mathematics and Information Science,}\\
{\small Minnan Normal University,}\\
{\small Zhangzhou, Fujian 363000, CHINA.}\\
{\small 2. Center for Discrete Mathematics and Theoretical Computer Science,}\\
{\small Fuzhou University,}
\\{\small Fuzhou, Fujian 350003, CHINA.}\\
{\small e-mail: laich@winmail.cn; laichunhui@mnnu.edu.cn
} \\
{\small (Chunhui Lai, Corresponding author) } \\
{\small liumingjing1983@163.com (Mingjing Liu) }}
\date{}
\maketitle
\begin{center}
\begin{minipage}{4.3in}
\vskip 0.1in
\begin{center}{\bf Abstract}\end{center}
 { Let $f(n)$ be the
maximum number of edges in a graph on $n$ vertices in which no two
cycles have the same length. Erd\"{o}s raised the problem of
determining $f(n)$. Erd\"{o}s conjectured that there exists a
positive constant $c$ such that $ex(n,C_{2k})\geq cn^{1+1/k}$.
Haj\'{o}s conjecture that every simple even graph on $n$ vertices can
be decomposed into at most $n/2$ cycles.  We present the problems,
conjectures related to these problems and we summarize the know results.
We do not think Haj\'{o}s conjecture is true. } \par
\vspace*{0.1in}
 \par
 \vspace*{0.1in}
   \par
\end{minipage}
\end{center}
\par
\section*{Erd\"{o}s Problem 1}\par
\par
Let $f(n)$ be the maximum number of edges in a graph on $n$ vertices
in which no two cycles have the same length. In 1975, Erd\"{o}s
raised the problem of determining $f(n)$ (see Bondy and Murty [1],
p.247, Problem 11).  Shi[41] proved a lower bound.
\par
\par
\noindent{\bf Theorem 1.1 (Shi[41])}
$$f(n)\geq n+
 [(\sqrt {8n-23} +1)/2]$$ for $n\geq 3$.
\par
 Chen, Lehel, Jacobson, and Shreve [4], Jia[25],
 Lai[28,29,30,31,32,33,34], Shi[42,43,44,45,46,47],
 Y. Shi, Y. Tang, H. Tang, L. Gong and L. Xu[48]
 obtained some results.
\par
 Boros, Caro, F\"uredi and Yuster[3]
 proved an upper bound.
\par
\noindent{\bf Theorem 1.2 ( Boros, Caro, F\"uredi and Yuster[3])}
  For $n$ sufficiently large, $$f(n) < n+1.98\sqrt{n}.$$
\par
Lai [35] improved the lower bound.
\par
\noindent{\bf Theorem 1.3 (Lai [35])}
  $$f(n)\geq n+ \sqrt {2.4} \sqrt {n}(1-o(1))$$ and proposed the following conjecture:
  \par
  \noindent{\bf Conjecture 1.4 (Lai [35])}
   $$\lim_{n \rightarrow \infty} {f(n)-n\over \sqrt n}=\sqrt {2.4}.$$
\par
It seems difficult to prove this conjecture. It would be nice to
prove one of the following weaker conjectures:
\par
\noindent{\bf Conjecture 1.5 (Lai[30])}
$$\liminf_{n \rightarrow \infty} {f(n)-n\over \sqrt n} \leq \sqrt {3}.$$

  \par
 \noindent{\bf Conjecture 1.6 (Lai[31])}
$$\liminf_{n \rightarrow \infty} {f(n)-n\over \sqrt n} \leq \sqrt {2.4}.$$
 \par
 Let $f_2(n)$ be the maximum number of edges in a $2$-connected
 graph on $n$ vertices in which no two cycles have the same length.
 \par
  Shi[44] proved that
  \par
\noindent{\bf Theorem 1.7 (Shi[44])} For every integer $n\geq3$,
$f_{2}(n)\leq n+[\frac{1}{2}(\sqrt{8n-15}-3)].$
\par
 Chen, Lehel, Jacobson, and Shreve [4]
 proved that
 \par
\noindent{\bf Theorem 1.8 (Chen, Lehel, Jacobson, and Shreve [4])}
$f_{2}(n)\geq n+\sqrt{n/2}-o(\sqrt{n})$
\par
  Boros, Caro, F\"uredi and Yuster [3] improved this lower bound significantly:
\par
\noindent{\bf Theorem 1.9 (Boros, Caro, F\"uredi and Yuster [3])}
$f_{2}(n)\geq n+ \sqrt{n}-O(n^{\frac{9}{20}})$.
\par
\noindent{\bf Corollary 1.10 (Boros, Caro, F\"uredi and Yuster [3])}
$$\sqrt{2}\geq \limsup \frac{f_{2}
(n)-n}{\sqrt{n}}\geq\liminf\frac{f_{2} (n)-n}{\sqrt{n}}\geq 1.$$
\par
 Boros, Caro, F\"uredi and Yuster [3] made the following conjecture:
\par
\noindent{\bf Conjecture 1.11 ( Boros, Caro, F\"uredi and Yuster [3])}\ \
$$\lim \frac{f_{2} (n)-n}{\sqrt{n}}=1.$$
\par
It is easy to see that Conjecture 1.11 implies the (difficult) upper
bound in the Erd\"{o}s-Turan Theorem [9,12](see Boros, Caro, F\"uredi
and Yuster [3]).
\par
Markstr\"om [27] raised the following problem:
 \par
\noindent{\bf Problem 1.12 (Markstr\"om [27])} Determine
 the maximum number of edges in a hamiltonian
graph on $n$ vertices with no repeated cycle lengths.
\par
Let  $g(n)$ denote the least number of edges of a graph which
contains a cycle of length $k$ for every $1\leq k\leq n$. Jia[25]
proved the following results:
 \par
 \noindent {\bf Theorem 1.13 (Jia[25])}
  \par When $n$ is sufficiently large, $$n+\log_{2}n-1\leq g(n)\leq
n+\frac{3}{2}\log_{2}n+1.$$
\par

\noindent{\bf Theorem 1.14 (Jia[25])}\ \
   \par For a sufficiently large positive integer $n$, $g(n)\leq n+\log_{2}n+\frac{3}{2}\log_{2}\log_{2}n+O(1)$
\par

\noindent{\bf Corollary 1.15 (Jia[25])}\  \par For $n$ sufficiently large,
$g(n)=n+\log_{2}n+O(\log_{2}\log_{2}n)$.
\par
Jia[25] made the following conjecture:
\par
\noindent{\bf
Conjecture 1.16 (Jia[25])}\ \ $$g(n)=n+\log_{2}n+O(1),$$ as
$n\rightarrow\infty$.

\par
The sequence $(c_1,c_2,\cdots,c_n)$  is the cycle length
distribution of a graph $G$ of order $n$£¬  where  $c_i$ is the
number of cycles of length $i$ in $G$. Let $f(a_1,a_2,\cdots,$
$a_n)$ denote the maximum possible number of edges in a graph which satisfies
$c_i\leq a_i$  where $a_i$  is a nonnegative integer. Shi posed the
problem of determining $f(a_1,a_2,\cdots,a_n)£¬$   which extended
the problem due to Erd\"{o}s, it is clearly that $f(n)=f(1,1,\cdots,1)$
(see Xu and Shi[52]).
 \par
 The lower bound
$f(0,0,2,\cdots,2)$  is given by Xu and Shi[52].
 \par
\noindent{\bf Theorem 1.17 (Xu and Shi[52])}  For  $n\geq 3$,
    \par
    $f(0,0,2,...,2)\geq n-1+[(\sqrt {11n-20})/2],$
     \par and the equality holds when $3\leq n \leq 10.$
\par
Given a graph $H$, what is the maximum number of edges of a graph
with $n$ vertices not containing $H$ as a subgraph? This number is
denoted $ex(n,H)$, and is known as the Turan number.
 \par
 We denote by  $m_i(n)$ the numbers of cycles of length $i$
 in the complete graph $K_n$ on $n$ vertices. Obviously,
 $$ex(n,C_{k})$$
 $$=f(0,0,m_3(n),
 \cdots,$$ $$m_{k-1}(n),0,m_{k+1}(n),
 \cdots,m_n(n))$$ $$=f(0,0,2^{\frac{n(n-1)}{2}}, \cdots,$$ $$2^{\frac{n(n-1)}{2}},
 0,2^{\frac{n(n-1)}{2}},
  \cdots, 2^{\frac{n(n-1)}{2}}).$$
Therefore, finding  $ex(n,C_k)$ is a special case of determining
$f(a_1,a_2,\cdots,a_n)$.
 
 \par
\section*{ Erd\"{o}s conjecture 2} \par

\par
 P. Erd\"{o}s
conjectured that there exists a positive constant $c$ such that
$ex(n,C_{2k})\geq cn^{1+1/k}$(see Erd\"{o}s[11]). Erd\"{o}s [8] posed the
problem of determining $ex(n,C_{4})$.
\par
Erd\"{o}s [10](published without proof) and Bondy and Simonovits [2] obtained that
     \par
\noindent{\bf Theorem 2.1 (Erd\"{o}s [10] and Bondy and Simonovits [2])}
$$ex(n,C_{2k})\leq ckn^{1+1/k}$$
\par
 Wenger [49] proved the following:
     \par
\noindent{\bf Theorem 2.2 (Wenger [49])} $$ex(n,C_{4})\geq
(\frac{n}{2})^{3/2},$$ $$ex(n,C_{6})\geq (\frac{n}{2})^{4/3},$$
$$ex(n,C_{10})\geq (\frac{n}{2})^{6/5}.$$
\par
F\"uredi[18] proved that
     \par
\noindent{\bf Theorem 2.3 (F\"uredi[18])} If $q$ is a power of 2, then
$$ex(q^{2} + q + 1,C_{4})= q(q+1)^{2}/2.$$
 \par
F\"uredi[19] also showed the following:
      \par
\noindent{\bf Theorem 2.4 (F\"uredi[19])} Let $G$ be a quadrilateral-free graph
with $e$ edges on $q^{2} + q + 1$ vertices, and suppose that $q \geq
15.$ Then $e \leq q(q+1)^{2}/2.$
     \par
\noindent{\bf Corollary 2.5 (F\"uredi[19])} If $q$ is a prime power greater
than 13, $n= q^{2} + q + 1.$ Then
$$ex(n,C_{4})= q(q+1)^{2}/2.$$
 \par
F\"uredi, Naor and Verstraete[20] proved that
     \par
\noindent{\bf Theorem 2.6 (F\"uredi, Naor and Verstraete[20])}
$$ex(n,C_{6})> 0.5338n^{4/3}$$ for  infinitely many $n$ and
$$ex(n,C_{6})< 0.6272n^{4/3}$$ if $n$ is sufficiently large.
 \par
 This refute the Erd\"{o}s-Simonovits conjecture in 1982 for hexagons(see[20]).
\par
The survey article on this Erd\"{o}s conjecture can be found in
Chung[5].
 \par
There are not good sufficient and necessary conditions of  when a graph on $n$ vertices
contains k cycle. For $k=n$, it is Hamiltonnian problem, the survey article on
 Hamiltonnian problem can be found in Gould [21,22].
 \par

  \section*{ Haj\'{o}s conjecture 3}
\par

  \par
An \emph{eulerian graph} is a graph (not necessarily connected) in
which each vertex has even degree. Let $G$ be an eulerian graph. A
\emph{circuit decomposition} of $G$ is a set of edge-disjoint
circuits $C_1,C_2,\cdots,C_t$ such that $E(G)=C_1\cup
C_2\cup\cdots\cup C_t$. It is well known that every eulerian graph
has a circuit decomposition. A natural question is to find the
smallest number $t$ such that $G$ has a circuit decomposition of $t$
circuits?  Such smallest number $t$ is called the \emph{circuit
decomposition } number of $G$, denoted by $cd(G)$.  For each edge
$xy\in E(G)$, let $m(xy)$ be the number of edges between $x$ and
$y$. The \emph{multiple number} of $G$ is defined by
$m(G)=\sum\limits_{uv\in E(G)}(m(uv)-1).$ A reduction of a graph
$G$ is a graph obtained from $G$ by recursively applying the
following operations:
\par
1. Remove the edges of circuit.
\par
2. Delete an isolated vertex (0-vertex).
\par
3. Delete a 2-vertex with two distinct neighbors and add a new edge joining its two neighbors.
\par
4. If $u$ is a 4-vertex with 4 distinct neighbors $\{ x,y,z,w\}$ such that
$xy \in E(G)$ and $zw \notin E(G)$, then delete $u$ and joint $x$ and $y$
with a new parallel edge, and add a new edge between $z$ and $w.$
\par
A reduction is proper if it is not the original graph (see [16]). The following
conjecture is due to Haj\'{o}s(see [36]).
\par

\noindent\textbf{Haj\'{o}s conjecture:}  $$cd(G)\leq
\frac{|V(G)|}{2}$$ for every simple eulerian graph $G$.
\par
Lovasz [36] proved the following:
\par
\par
\noindent\textbf{Theorem 3.1 (Lovasz [36])} A graph of $n$ vertices
can be covered by $ \leq  [n/2]$ disjoint paths and circuits.

\par
Jiang [26] and Seyffarth [40] considered planar
eulerian graphs.
\par
\noindent\textbf{Theorem 3.2 (Jiang [26] and Seyffarth [40])}
  \emph{$cd(G)\leq \frac{|V(G)|-1}{2}$ for every simple planar
eulerian graph $G$.}
\par

\par
Granville and  Moisiadis [23] and  Favaron and Kouider [17] extended to multigraphs.
\par
\par
\noindent\textbf{Theorem 3.3 (Granville and Moisiadis [23] and
Favaron and  Kouider [17])}  If $G$ is an even multigraph of order
$n$, of size $m$, with $\Delta (G) \leq 4$, then \emph{$cd(G)\leq
\frac{n+M-1}{2}$   where $M=m-m^*$ and $m^*$ is the size of the
simple graph induced by $G$.}
\par
  Fan and  Xu [16] proved that:
\par
\noindent\textbf{Theorem 3.4 (Fan and Xu [16])} \emph{If $G$ is an
eulerian graph with $$cd(G)> \frac{|V(G)|+m(G)-1}{2}$$ then $G$ has
a reduction $H$ such that
$$cd(H)> \frac{|V(H)|+m(H)-1}{2}$$ and the number of vertices of
degree less than six in $H$ plus $m(H)$ is at most one.}
\par
\noindent\textbf{Corollary 3.5 (Fan and Xu [16])} \emph{Haj\'{o}s
conjecture is valid for projective graphs.}\par

\noindent\textbf{Corollary 3.6 (Fan and Xu [16])} \emph{Haj\'{o}s
conjecture is valid for $K_6^-$ minor free graphs.}
\par
 Xu[50] also proved the following two results:
 \par
\noindent\textbf{Theorem 3.7 (Xu[50])} \emph{If $G$ is an eulerian
graph with $$cd(G)> \frac{|V(G)|+m(G)-1}{2}$$ such that
$$cd(H)\leq \frac{|V(H)|+m(H)-1}{2}$$  for each proper reduction of
$G$, then $G$ is 3-connected. Moreover, if $S=\{x,y,z\}$  is a 3-cut
of $G$, letting $G_1$ and $G_2$ be the two induced subgraph of $G$
such that
 $V(G_1)\cap V(G_2)=S$ and $E(G_1)\cup E(G_2)=E(G)$,   then either $S$
 is not an independent set, or $G_1$ and $G_2$ are both eulerian graphs.}
 \par
 \noindent\textbf{Corollary 3.8 (Xu[50])} \emph{ To prove Haj\'{o}s' conjecture, it suffices to prove
 $$cd(G)\leq
\frac{|V(G)|+m(G)-1}{2}$$  for every 3-connected  eulerian graph
$G$.}
\par
 Fan [14] proved that
\par
\noindent\textbf{Theorem 3.9 (Fan [14])} Every eulerian graph on $n$
vertices can be covered by at most $\lfloor \frac{n-1}{2}\rfloor $
circuits such that each edge is covered an odd number of times.
\par
This settles a conjecture made by Chung in 1980(see[14]).
\par
Xu and Wang[51] gave the following result:
\par
\noindent\textbf{Theorem 3.10 (Xu and Wang[51])} \emph{The edge set
of each even toroidal graph can be decomposed into at most $(n+3)/2$
circuits in $O(mn)$ time, where a toroidal graph is a graph
embedable on the torus.}
\par
\noindent\textbf{Theorem 3.11 (Xu and Wang[51])} \emph{The edge set
of each  toroidal graph can be decomposed into at most $3(n-1)/2$
circuits and edges in $O(mn)$ time.}
\par
\noindent{\bf We do not think Haj\'{o}s conjecture is true.}
\par
By the proof of Lemma 3.3 in N. Dean [What is the smallest number of dicycles in a
dicycle decomposition of an Eulerian digraph?
J. Graph Theory 10 (1986), no. 3, 299--308.],
if exists $k$ vertices counterexample,
 then will exist $ik-i+1(i=2,3,4,...)$ vertices counterexamples.
\par
A related problem was conjectured by  Gallai (see [36]):
\par
\noindent{\bf Conjecture 3.12 (Gallai's conjecture)}
Every simple connected graph on $n$ vertices can be
decomposed into at most $(n+1)/2$ paths.
\par
Lovasz [36] proved that
\par
\par
\noindent\textbf{Theorem 3.13 (Lovasz [36])} If a graph has $u$ odd
vertices and $g$ even vertices ($g \geq 1)$, then it can be covered
by $u/2 + g - 1$ disjoint paths.
\par
\noindent\textbf{Theorem 3.14 (Lovasz [36])} Let a locally finite graph
 have only odd vertices. Then it  can be covered by disjoint finite paths
 so that every vertex is the endpoint of just one covering path.
\par
The path number of a graph $G$, denoted $p(G)$, is the minimum number of
edge-disjoint paths covering the edges of $G$. Donald [7]
proved the following:
\par
\noindent\textbf{Theorem 3.15 (Donald [7])} If a graph with $u$
vertices of odd degree and $g$ nonisolated vertices of even degree.
Then
 $$p(G) \leq u/2 + [\frac{3}{4}g] \leq [\frac{3}{4}n].$$
\par
\par
Pyber [37] proved that
\par
\par
\noindent\textbf{Theorem 3.16 (Pyber [37])} A graph $G$ of  $n$
vertices can be covered by $n-1$ circuits and edges.
\par
\noindent\textbf{Theorem 3.17 (Pyber [37])} Let $G$ be a graph of $n$ vertices
and $\{C_1,\ldots, C_k \}$ be a set of circuits and edges such that
$\cup^{k}_{i=1}E(C_i)=E(G)$ and $k$ is minimal. Then we can choose $k$ different edges,
 $e_i \in E(C_i),$ such that these edges form a forest in $G$.
\par
\noindent\textbf{Theorem 3.18 (Pyber [37])} Let $G$ be a graph of $n$ vertices
not containing $C_4$. Then $G$ can be covered by $[(n-1)/2]$ circuits and $n-1$ edges.
\par
Pyber [38] proved that
\par
\noindent\textbf{Theorem 3.19 (Pyber [38])}  Every connected  graph $G$ on $n$
vertices can be covered by $n/2 + O(n^{3/4})$ paths.
\par
\par
\noindent\textbf{Theorem 3.20 (Pyber [38])}  Every connected  graph on $n$
vertices with $e$ edges can be covered by $n/2 + 4(e/n)$ paths.
\par
\par
Reed [39]
proved that
\par
\noindent\textbf{Theorem 3.21 (Reed [39])}  Any connected cubic graph $G$ of
order $n$ can be covered by $\lceil n/9\rceil$ vertex disjoint paths.
\par
\par
Dean, Kouider [6] and Yan [53] proved that
\par
\noindent\textbf{Theorem 3.22 (Dean, Kouider [6] and Yan [53])} If a graph(possibly disconnected) with $u$
vertices of odd degree and $g$ nonisolated vertices of even degree.
Then
 $$p(G) \leq u/2 + [\frac{2}{3}g].$$
\par
\par
Fan [13] proved that
\par
\noindent\textbf{Theorem 3.23 (Fan [13])}  Every connected  graph on $n$
vertices can be covered by at most $\lceil n/2 \rceil$ paths.
\par
This settles a conjecture made by Chung in 1980(see[13]).
\par
\par
\noindent\textbf{Theorem 3.24 (Fan [13])}  Every 2-connected  graph on $n$
vertices can be covered by at most $\lfloor \frac{2n-1}{3} \rfloor$ circuits.
\par
This settles a conjecture made by Bondy in 1990(see[13]).
\par
\par
\noindent\textbf{Corollary 3.25 (Fan [13])}  Let $G$ be a  2-edge-connected  graph on $n$
vertices. Then $G$  can be covered by at most $\lfloor \frac{3(n-1)}{4} \rfloor$ circuits.
\par
\par
Fan [15] define a graph operation,
called $\alpha$-operation and proved that
\par
\noindent\textbf{ Definition 3.26(Fan [15])} Let H be a graph. A pair $(S, y)$, consisting of an independent set $S$ and a
vertex $y \in S,$ is called an $\alpha$-pair if the following holds: for every vertex $v \in S \setminus \{y\},$ if
$d_H (v) \geq 2,$ then (a) $d_H (u) \leq 3$ for all $u \in N_H (v)$ and (b) $d_H (u) = 3$ for at most two vertices
$u \in N_H (v).$ (That is, all the neighbors of $v$ have degree at most 3, at most two of which have
degree exactly 3.) An $\alpha$-operation on $H$ is either (i) add an isolated vertex or (ii) pick an
$\alpha$-pair $(S, y)$ and add a vertex $x$ joined to each vertex of $S$, in which case the ordered triple
$(x, S, y)$ is called the $\alpha$-triple of the $\alpha$-operation.
\par
\noindent\textbf{ Definition 3.27(Fan [15])} An $\alpha$-graph is a graph that can be obtained from the empty set via a sequence
of $\alpha$-operations.
\par
\noindent\textbf{Theorem 3.28 (Fan [15])}  Let $G$ be a  graph on $n$
vertices ( not necessarily connected). The E-subgraph of $G$ is the subgraph induced by
 the vertices of even degree in $G$. If the E-subgraph of $G$
is an $\alpha$-graph, then $G$ can be decomposed into  $\lfloor n/2 \rfloor$ paths.
\par
\par
\noindent\textbf{Corollary 3.29 (Fan [15])}  Let $G$ be a  graph on $n$
vertices ( not necessarily connected). If each block of the E-subgraph of $G$
 is a triangle-free graph with maximum degree at most 3, then $G$ can be
 decomposed into  $\lfloor n/2 \rfloor$ paths.

\par
Harding and McGuinness [24] proved that
\par
\par
\noindent\textbf{Theorem 3.30 (Harding and McGuinness [24])} For every simple graph $G$ have girth $g \geq 4$, with $u$
vertices of odd degree and $w$ nonisolated vertices of even degree, there is a path-decomposition having at most
$u/2 + \lfloor \frac{g+1}{2g}w \rfloor $ paths.
\par

\par
\section*{Acknowledgment}
The authors would like to thank Professor Y. Caro, G. Fan, R. Gould,
 B. Xu, R. Yuster for their advice and sending some papers to us.
The authors would like to thank the referee for his many valuable comments and suggestions.
\par

\end{document}